\documentclass[a4paper,12pt]{article}
  \usepackage{macros}

\title{Sur les remplissages holomorphes équivariants}
\author{Beno\^{\i}t Kloeckner}

\newcommand{\SSD}{{\mathcal S}}
\newcommand{\BS}{{\mathcal B}}
\newcommand{\CBS}{\CP^2\priv{\mathcal B}}

\begin{document}

\maketitle

\section{Introduction}

\subsection{Remplissages équivariants}

La notion de remplissage d'une structure géométrique est très largement
étudiée. Nous nous proposons d'en étudier une version dynamique dans
un cadre CR. Introduisons pour cela une définition.
\begin{defi}[remplissage équivariant]
Soit $M$ une variété CR de dimension $2n-1$ et $F$ 
un sous-groupe de $\Aut(M)$.
On appelle \defini{remplissage holomorphe}
(ou simplement \defini{remplissage}) de $M$ toute variété complexe
à bord $R$ de dimension $n$ telle qu'il existe
un isomorphisme CR 
\[d:M\rightarrow \bord R\]

En conjuguant par $d$, l'action de $F$ sur $M$ donne une action
de $F$ sur $\bord R$. Si elle se prolonge en une action sur $R$ par 
biholomorphismes, on dit que $R$ est un remplissage de $M$
\defini{équivariant relativement à} $F$.
\end{defi}
Comme on ne précise pas l'isomorphisme $d$ dans un remplissage, ce n'est
en fait que la classe de conjugaison de $F$ dans $\Aut(M)$ qui nous importe.
Une définition précise d'une variété complexe à bord est donnée à la section
\ref{varietesAbord}.

Étant donnée une variété CR sur laquelle agit un groupe, on peut
chercher à classer les remplissages de la variété qui sont équivariants
relativement à l'action. Dans cet article, on se restreint au
cas des variétés strictement pseudoconvexes.
Parmi celles-ci, seule la sphère standard admet une dynamique riche 
comme le montre le théorème suivant. On appelle sphère standard et 
on note $\SSD^{2n-1}$ la variété
CR abstraite isomorphe à la sphère unité
$\{\sum_{i=1}^n \abs{z_i}^2=\abs{z_0}\}$
de $\CP^n$, muni de coordonnées homogènes $[z_0:\dots:z_n]$.
\begin{theo}[Webster--Schoen \cite{Schoen}]
Soit $M$ une variété CR compacte et strictement pseudoconvexe de dimension
$2n-1$. Si $\Aut(M)$ est non compact, alors $M$ est CR-équivalente
à la sphère standard $\SSD^{2n-1}$.
\end{theo}

On se concentre par conséquent sur la sphère standard et 
les sous-groupes fermés non compacts de son groupe d'automorphismes
$\SU{1,n}$.
De plus on se restreint au cas $n=2$.

Comme la sphère standard est strictement pseudoconvexe, ses deux 
côtés ne jouent pas des rôles symétriques. Un remplissage $R$
est dit \defini{convexe} ou \defini{concave} selon le côté
du bord occupé par l'intérieur de $R$.

Plongée comme sphère unité dans $\CP^2$,
$\SSD^{2n-1}$ y découpe deux composantes connexes non isomorphes. On note 
$\BS$ la boule standard définie par l'équation
$\{\sum_{i=1}^n \abs{z_i}^2<\abs{z_0}\}$
et $\CBS$ son complémentaire. Ces deux remplissages sont respectivement
convexe et concave. Ils sont tous deux équivariants relativement à $\SU{1,2}$
car l'action de celui-ci sur la sphère
se prolonge à tout $\CP^2$.

Le cas des remplissages convexes équivariant est assez simple.
\begin{prop}
Soit $R$ un remplissage convexe de $\SSD^3$ équivariant relativement
à un sous-groupe fermé non compact de $\SU{1,2}$.
Alors $R$ est isomorphe à la boule standard $\adherence{\BS}$.
\end{prop}

Notons que sans l'hypothèse d'équivariance, on dispose de résultats 
importants (et beaucoup plus profonds). Citons en particulier
une conséquence des remplissages par disques holomorphes.

\begin{theo}[Eliashberg \cite{Eliashberg}]
Soit $M$ une variété CR de dimension $3$ difféomorphe à la sphère
et strictement pseudoconvexe.
Un remplissage convexe de $M$ est difféomorphe à la boule fermée
éventuellement éclatée en quelques points.
\end{theo}

Sous l'hypothèse d'équivariance, le cas concave est un peu plus 
subtil que le cas convexe. On ne peut plus espérer
l'unicité car on construit facilement différents exemples à
partir de $\CBS$. Choisissons un élément hyperbolique ou parabolique
$\phi$ de $\SU{1,2}$ (c'est-à-dire un élément qui engendre un sous-groupe
fermé non compact). Alors $\phi$ agit sur $\CBS$ et s'il admet un point
fixe $p$ en dehors du bord, on peut éclater $\CBS$ en $p$ pour obtenir
un nouvel exemple. En effet $\phi$ agit encore sur l'éclaté et engendre
un groupe fermé non compact. On peut continuer à éclater des points fixes
de l'action de $\phi$ pour obtenir de nombreux exemples. Ces remplissages
concaves de $\SSD^3$ sont appelés \defini{éclatements standards} de $\CBS$.

Notre but est de démontrer le résultat suivant.
\begin{theo}\label{theoremplissages}
Considérons la sphère standard $\SSD^3$ et
un sous-groupe $F$ de $\Aut(\SSD^3)=\SU{1,2}$.
Si $F$ est non compact et fermé dans $\Aut(\SSD^3)$, alors tout 
remplissage concave de $\SSD^3$ 
équivariant relativement à $F$ est un éclatement standard de $\CBS$.
\end{theo}

Notons qu'il est vain d'espérer classifier précisément
les variétés complexes à bord strictement pseudoconcave
sans l'hypothèse d'équivariance. 
En effet priver n'importe quelle variété d'un domaine strictement 
pseudoconvexe inclu dans l'une de ses cartes fournit une grande 
diversité d'exemples.

Pour se convaincre que demander seulement la non-compacité du groupe 
d'automorphismes holomorphes du remplissage ne rend pas le problème
accessible, considérons un exemple édifiant
dû à McMullen.

\subsection{L'exemple de McMullen}

Dans \cite{McMullen}, McMullen 
décrit une surface K3 (et plus récemment dans \cite{McMullen2},
des surfaces rationnelles) possédant un
automorphisme $\phi$ cumulant les deux propriétés suivantes :
\begin{enumerate}
\item il est d'entropie positive ;
\item il admet un domaine de Siegel.
\end{enumerate}
Comme $\phi$ est d'entropie positive, il engendre un groupe non compact.
Mais le domaine de Siegel est par définition un ouvert de la surface
dans lequel $\phi$ est conjugué à une rotation, donc il contient
un ouvert stable biholomorphe à la boule standard. Si on prive la surface
de cet ouvert, on obtient une surface à bord strictement pseudoconcave
possédant un automorphisme d'entropie positive. L'image du groupe
engendré par $\phi$ dans le groupe d'automorphismes du bord est
relativement compacte et l'inclusion ne peut pas être propre.

En un sens, on peut dire que la dynamique de la surface ainsi construite
reste à l'écart du bord.

\subsection{Reformulations} 

Dans cette section, nous proposons une seconde interprétation
du théorème \ref{theoremplissages}.

Remarquons tout d'abord que le groupe des automorphismes
d'une variété complexe à bord s'injecte naturellement
dans le groupe des automorphismes CR du bord.
Toutefois cette injection, comme on le verra plus bas, n'est pas
nécessairement propre.

\begin{coro}
Soit $X$ une surface complexe à bord strictement pseudoconcave.
Si 
\begin{enumerate}
\item le groupe d'automorphisme de $X$ est non compact et 
\item l'inclusion de $\Aut(X)$ dans $\Aut(\bord X)$
      est propre 
\end{enumerate}
alors $X$ est un éclaté de $\CP^2\priv\BS$, obtenu
comme dans le théorème \ref{theoremplissages}.
\end{coro}

Cet énoncé découle du théorème \ref{theoremplissages} grâce
au théorème de Webster--Schoen : nos hypothèses assurent que
le groupe d'automorphismes du bord est non compact, donc que
le bord est bien la sphère standard.

On peut encore reformuler le résultat en terme de sphère invariante.

\begin{coro}
Soit $X$ une surface complexe sans bord.
Si
\begin{enumerate}
\item  $X$ possède une hypersurface réelle plongée $H$  qui la déconnecte et
\item le groupe des automorphismes de $X$ qui préservent $H$ y induit
      un groupe fermé non compact d'automorphismes CR
\end{enumerate}
alors $X$ est nécessairement un éclaté de $\CP^2$, obtenu
comme dans le théorème \ref{theoremplissages}.
\end{coro}

Pour se ramener à cette forme,
il suffit de constater que le remplissage concave peut être
complété par un remplissage convexe par la boule pour donner une variété
complexe sans bord (voir la démonstration plus loin).

\subsection{Variétés complexes à bord}\label{varietesAbord}

On peut généraliser de la façon suivante la notion de variété
complexe.

\begin{defi}[variété complexe à bord]
Soient $X$ un espace topologique paracompact séparé et $n$ un entier positif.
On appelle \defini{carte complexe} (de dimension $n$) de $X$
tout quadruplet ${\mathscr U}=(U,\varphi,H,A)$ où
\begin{enumerate}
\item $U$ est un ouvert de $X$, appelé \defini{domaine} de la carte ;
\item $A$ est un ouvert de $\mC^n$ et :\begin{itemize}
  \item soit $H$ est une hypersurface analytique 
      réelle de $A$ et $A\priv H$ a deux composantes connexes $D_+$ et 
      $D_-$, 
  \item soit $H$ est vide et on note $D_+=A$ ;\end{itemize}
\item $\varphi$ est un homéomorphisme de $U$ vers $D_+\cup H$.
\end{enumerate}
Deux cartes complexes $(U,\varphi,H,D)$ et $(V,\psi,K,E)$ sont dites
\defini{compatibles} si le changement de coordonnées
\[\psi\comp\varphi^{-1} : \varphi(U\cap V)  \longrightarrow \psi(U\cap V)\]
se prolonge en un biholomorphisme dans un voisinage de son domaine
de définition vers un voisinage de son image. On considère
cette condition satisfaite si $U\cap V=\varnothing$.
Un \defini{atlas complexe} est un ensemble de cartes compatibles
dont les domaines recouvrent $X$.

On appelle \defini{variété complexe à bord}
\index{Variete complexe a bord@Variété complexe à bord}
de dimension $n$ un espace topologique paracompact $X$
muni d'un atlas complexe.
\end{defi}

Enrichir un atlas d'une carte compatible à toutes les précédentes
est une opération anodine ; on dit que deux atlas sont équivalents
si leur union est encore un atlas et on identifie en fait
la structure de variété complexe à bord à la classe d'équivalence
de l'atlas choisi.

Une variété complexe à bord est naturellement munie d'une structure
de variété différentielle à bord.

L'intérieur d'une variété complexe à bord est muni d'une
structure complexe tandis que son bord porte une
structure CR, définie localement par les cartes.

D'après le théorème de Newlander--Nirenberg, une variété réelle
$X$ munie d'un opérateur complexe ($J:TX\rightarrow TX$, $J^2=-\id$)
défini et intégrable jusqu'à un voisinage du bord est une variété complexe
à bord.

\begin{defi}[isomorphisme]
Soient $X$ et $Y$ deux variétés complexes à bord. Un homéomorphisme
$F:X\rightarrow Y$ est un \defini{isomorphisme complexe}
\index{Isomorphisme complexe} si pour toute carte ${\mathscr V}$
de $Y$, $F\pullback {\mathscr V}$ est compatible
à toute carte de $X$.
Si $X=Y$, on parle d'\defini{automorphisme}\index{Automorphisme complexe}.
\end{defi}

Un isomorphisme n'est donc rien d'autre qu'un biholomorphisme qui se 
prolonge à un voisinage du bord.

\section{Dynamique de $\SU{1,2}$}

Nous aurons besoin de bien comprendre les éléments hyperboliques et 
paraboliques de $\SU{1,2}$ pour mener à bien la démonstration du
théorème \ref{theoremplissages}. En particulier, nous sommes intéressés
par leurs bassins d'attraction et de répulsion.

Dans la suite on utilise sur $\CP^2$ des coordonnées homogènes
$[x:y:z]$ et on note $Q$ la forme hermitienne 
$\abs{x}^2+\abs{y}^2-\abs{z}^2$. La boule unité $\BS$ s'écrit ainsi
$\{Q<0\}$.

Le théorème du point fixe de Brouwer implique que chaque élément de
$\SU{1,2}$ admet un point fixe dans la boule fermée $\adherence{\BS}$.
On peut les classer en trois
catégories : un élément de $\SU{1,2}$ est 
\begin{itemize}
\item \defini{elliptique}\index{Isométrie!elliptique}
      s'il fixe un point de $\BS$ ;
\item \defini{parabolique} s'il ne fixe aucun point de $\BS$ mais 
      exactement un point de $\bord \BS$ ;
\item \defini{hyperbolique} s'il ne fixe aucun point de $\BS$ mais
      exactement deux points de $\bord \BS$,
\end{itemize}
et tout élément est d'un de ces trois types.

Remarquons qu'un point de $\CP^2$ fixé par $A\in\SU{1,2}$ correspond
à un vecteur propre de la matrice $A$.

\subsection{Description de l'algèbre de Lie}

L'algèbre de Lie contient beaucoup d'information sur le groupe,
et se prête souvent mieux aux calculs.
Les éléments de $\su{1,2}$ sont les matrices $3\times3$ de
la forme
\begin{equation}
A=
\begin{pmatrix}
-i(b_1+b_2)    &    l_1        & l_2 \\
\conjugue{l_1} &    ib_1       & c  \\
\conjugue{l_2} & -\conjugue{c} & ib_2 
\end{pmatrix}
\label{LieAlgebraForm}
\end{equation}
où $b_1$ et $b_2$ sont des nombres réels tandis que $l_1$, $l_2$ et
$c$ sont des nombres complexes.


Un élément de $\su{1,2}$ est dit elliptique, parabolique ou
respectivement hyperbolique si l'adjectif s'applique à son exponentielle.

Nous nous intéressons maintenant aux transformations paraboliques
et hyperboliques, qui nous seront utiles par la suite.

\subsection{Éléments hyperboliques}

Par définition, un élément hyperbolique $A\in\SU{1,2}$ a deux points
fixes $p_1$ et $p_2$ dans $\bord \BS$. Par chacun de ses points, il
passe exactement une droite projective tangente à la sphère $\bord \BS$,
dont la direction est la droite complexe $\xi_{p_i}\subset T_{p_i}\bord \BS$ 
de la distribution de contact sous-jacente à la structure CR.
Ces deux droites $L_1$, $L_2$ se coupent en un point $q$ de $\CP^2$ et comme
elles sont globalement préservées, $q$ est un point fixe pour l'action
projective de $A$. Ces trois points non alignés correspondent
à trois vecteurs propres linéairement indépendants, il n'y en a donc
pas d'autre.

\begin{figure}[htb]\begin{center}
\includegraphics{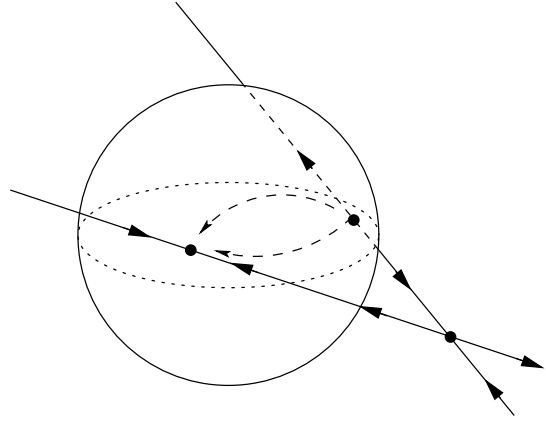}
\caption{Dynamique d'une transformation 
           hyperbolique}\label{trans_hyper}
\end{center}\end{figure}

L'action de $\SU{1,2}$ est deux fois transitive sur $\bord \BS$ ; quitte
à conjuguer $A$ on peut donc supposer que les $p_i$ sont les points
$[1:1:0]$ et $[1:-1:0]$. Alors $q$ doit être le point $[0:0:1]$,
intersection des droites $\{x=y\}$ et $\{x=-y\}$. En utilisant 
\refeq{LieAlgebraForm} on peut écrire $A$ sous la forme
\begin{equation}
A=\exp\begin{pmatrix}
ib & l  & 0   \\
l  & ib & 0   \\
0  & 0  & -2ib
\end{pmatrix}
\end{equation}
où $l\neq 0$ et $b$ sont des nombres réels. Quitte à considérer $-A$
qui a la même action sur $\CP^2$, on peut supposer $l>0$.
Les valeurs propres de $A$ sont $\exp(l+ib)$, $\exp(-l+ib)$,
$\exp(-2ib)$ et $A$ est diagonalisable. Un calcul en coordonnées
projectives permet de déterminer la nature des trois points fixes.
L'un des points fixes du bord, disons $p_1$, est attractif avec
valeurs propres $\exp(-2l)$ et $\exp(-l-3ib)$. L'autre est répulsif
avec valeurs propres $\exp(2l)$ et $\exp(l-3ib)$. Enfin, $q$ est
hyperbolique avec valeurs propres $\exp(-l+3ib)$ et $\exp(l+3ib)$.
Sur la figure \ref{trans_hyper}, la sphère
           $\bord \BS$ est représentée en dimension $2$ au lieu de $3$ 
           et les droites complexes sont représentées par des 
           lignes.

Le bassin d'attraction de $p_1$ est $\CP^2\priv L_2$, le bassin
de répulsion de $p_2$ est $\CP^2\priv L_1$. L'union des deux bassins
est donc $\CP^2\priv\{q\}$.

\subsection{Éléments paraboliques}

Considérons maintenant un élément parabolique $a\in\su{1,2}$ et
$A$ son exponentielle.
Quitte à conjuguer, on peut supposer que son unique point fixe au bord 
est $p=[1:1:0]$.

L'action de $A$ sur la droite projective $L$ passant
par $p$ et tangente à $\bord \BS$ ne peut être hyperbolique, donc la valeur
propre de $a$ associée au vecteur propre $(1,1,0)$ est imaginaire pure
et on peut écrire
\begin{equation}
a=\begin{pmatrix}
-i(d_1+d_2)                      & i\left(d_1+\frac{1}{2}d_2\right) & c \\
-i\left(d_1+\frac{1}{2}d_2\right) & id_1                            & c \\
\conjugue{c}                     & -\conjugue{c}                   & id_2
\end{pmatrix}
\end{equation}
où $d_1$ et $d_2$ sont réels et $c$ est complexe.
Les valeurs propres de $a$ sont $id_2$ et $-id_2/2$ avec multiplicité $2$.
On peut maintenant distinguer trois cas, selon la forme de Jordan de
$A$.

\begin{figure}[htb]\begin{center}
\includegraphics{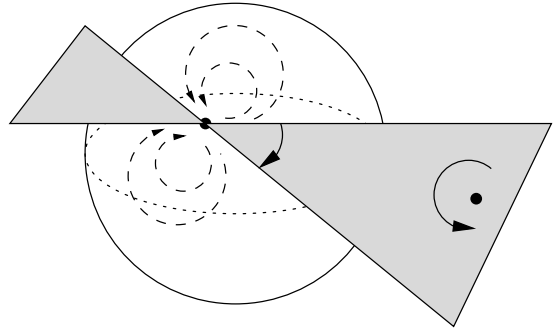}
\caption{Dynamique d'une transformation 
           parabolique à deux points fixes}\label{trans_para1}
\end{center}\end{figure}

Si $d_2\neq 0$, $p$ définit l'unique direction propre associée à la valeur
propre $-id_2/2$ donc $a$ n'est pas diagonalisable et la forme
de Jordan de $A$ est
\begin{equation}
\begin{pmatrix}
e^{-i\frac{d_2}{2}}  &     1           &  0 \\
0                & e^{-i\frac{d_2}{2}} &  0 \\
0                &     0           & e^{id_2}
\end{pmatrix}
\end{equation}
L'action de $A$ sur $L$ est alors conjuguée à
une rotation de $\CP^1$ et $A$ admet
un deuxième point fixe $q\in L$.
Sur la figure \ref{trans_para1}, on a cette fois représenté $L$
par une section plane pour y montrer la  rotation.

Dans ce cas le bassin d'attraction de $p$ est $\CP^2\priv L$.

\begin{figure}[htb]\begin{center}
  \includegraphics{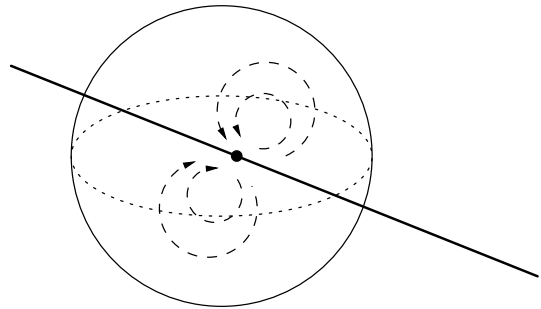}
  \caption{Dynamique d'une transformation parabolique fixant chaque
           point d'une droite}\label{trans_para2}
\end{center}\end{figure}

Si $d_2=0$ et $c=0$, $a$ est nilpotente d'ordre $2$, la forme
de Jordan de $A$ est
\begin{equation}
\begin{pmatrix}
1 & 1 & 0 \\
0 & 1 & 0 \\
0 & 0 & 1
\end{pmatrix}
\end{equation}
et $A$ fixe chaque point de $L$ (voir la figure \ref{trans_para2}).

Dans ce cas le bassin d'attraction de $p$ est encore $\CP^2\priv L$.

\begin{figure}[htb]\begin{center}
  \includegraphics{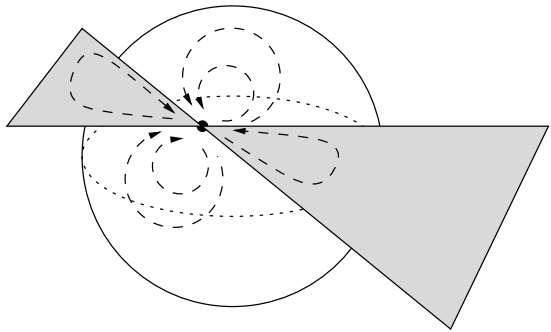}
  \caption{Dynamique d'une transformation parabolique fixant un seul
           point}\label{trans_para3}
\end{center}\end{figure}

Enfin si $d_2=0$ et $c\neq0$, $a$ est nilpotente d'ordre $3$, la
forme de Jordan de $A$ est
\begin{equation}
\begin{pmatrix}
1 & 1 & 0 \\
0 & 1 & 1 \\
0 & 0 & 1
\end{pmatrix}
\end{equation}
et $p$ est son seul point fixe, à la fois attractif et répulsif
(voir la figure \ref{trans_para3}).

Dans ce cas le bassin d'attraction est tout $\CP^2$.

On obtient donc le résultat suivant qui est une des clés
du théorème \ref{theoremplissages}.
\begin{lemm}\label{bassins}
Soit $g$ un élément parabolique ou hyperbolique de $\SU{1,2}$.
Notons $p_+$ le point fixe attractif de $g$,
$p_-$ son point fixe répulsif s'il est différent et $L_+$
la droite projective tangente à la sphère unité en $p_+$. 

Tout ouvert de $\CP^2$ contenant l'union du bassin 
d'attraction de $p_+$ et du bassin de
répulsion de $p_-$ contient :
\begin{enumerate}
\item un voisinage de la boule unité fermée et
\item la famille des droites projectives passant par $p_+$, sauf
      éventuellement $L_+$.
\end{enumerate}
\end{lemm}

\section{Remplissages de la sphère}

Dans cette dernière section nous donnons une démonstration du
théorème \ref{theoremplissages}. Commençons par en présenter
le déroulement.

On considère un remplissage $R$ de la sphère standard,
équivariant relativement à un sous-groupe $F\subset\SU{1,2}$
fermé et non compact. On note $d$ un isomorphisme
$\SSD\rightarrow \bord R$.

Notre premier outil est le théorème \ref{casferme} démontré
en annexe, qui implique que $F$ n'est pas 
purement elliptique. Il contient donc au moins une transformation 
hyperbolique ou parabolique, qui agit sur $R$ et sur $\CP^2$ après
conjugaison par $d^{-1} : \bord R\rightarrow \SSD\subset\CP^2$. 
On note $\phi$ l'automorphisme de $R$ et 
$\wt{\phi}$ celui de $\CP^2$.

On se base alors sur l'étude des bassins d'attraction et de répulsion
de $\wt{\phi}$. Elle permet de montrer de façon très élémentaire
que $R$ se complète par l'ajout d'une boule en une surface sans bord $X$
sur laquelle $\phi$ agit encore.
Mieux, elle permet d'exhiber dans $X$ une famille de
courbes rationnelles d'auto-intersection $1$, homologues entre elles
et dont la classe est invariante sous l'action de $\phi$.
Il découle alors du lemme de Noether que $X$ est une surface 
rationnelle. 

On utilise ensuite le théorème d'indice de Hodge pour montrer
que quitte à élever $\phi$ à une certaine puissance $\phi^k$,
on peut supposer que $\phi$ préserve globalement chaque courbe
exceptionnelle, ce qui permet de les contracter.

Il suffit alors d'étudier le cas des surfaces rationnelles minimales.

\subsection{Prolongement convexe}

\subsubsection{Prolongement des conjugaisons}

On dit de deux automorphismes $\psi$ et $\wt{\psi}$ de 
variétés complexes (éventuellement à bord) $X$ et $Y$
qu'ils sont \defini{localement conjugués} s'il existe
des ouverts $U\subset X$ et $V\subset Y$ et un biholomorphisme
$F:U\rightarrow V$ tel que $F\comp\psi=\wt{\psi}\comp F$ là où
cette expression est définie.
On dit que $U$ et $V$ sont les \defini{ouverts de conjugaison}.

On va utiliser notre étude des bassins d'attraction grâce
au lemme élémentaire qui suit.
\begin{lemm}[de prolongement]
Soient $\psi$ et $\wt{\psi}$ des automorphismes de variétés complexes
(éventuellement à bord) $X$ et $Y$, localement conjugués dans des ouverts
$U\subset X$ et $V\subset Y$. On peut alors prolonger la conjugaison
à des ouverts stables $U'\supset U$ et $V'\supset V$ où $V'$ contient
tous les bassins d'attraction et de répulsion des points fixes
de $\wt{\psi}$ dans $V$.
\end{lemm}

\preuve
notons $U_0=U$, $V_0=V$ et $F_0=F$. On définit récursivement
\begin{eqnarray}
U_{k+1}&=&\psi^{-1}(U_k)\nonumber\\
V_{k+1}&=&\wt{\psi}^{-1}(V_k)\nonumber\\
F_{k+1}&=&\wt{\psi}^{-1}\comp F_k\comp\psi
\nonumber
\end{eqnarray}
Ainsi  $F_k$ est une conjugaison locale entre les ouverts $U_k$ et $V_k$.
De plus là où deux $F_k$ sont définis simultanément, ils coïncident
par équivariance.
En passant à l'union, on peut donc construire une conjugaison
$F_\infty$ entre les ouverts
$\bigcup_k U_k$ et $\bigcup_k V_k$.
Ce dernier contient tous les bassins d'attractions des point fixes
de $\wt{\psi}$ contenus dans $V$. Comme $F$ est également une 
conjugaison de $\psi^{-1}$ et $\wt{\psi}^{-1}$ on peut à nouveau la
prolonger jusqu'à englober les bassins de répulsion. Par construction
les ouverts de conjugaison sont stables sous l'action de $\psi$
et $\wt{\psi}$ respectivement. 
\finpreuve

\subsubsection{Conséquences}

Le plongement de $\SSD^3$ (avec sa structure CR standard)
dans une surface est localement unique d'après le résultat suivant. 
\begin{theo}[Pinchuk \cite{SCV}]\label{pin}
Une application CR, \diffb{1} et non cons\-tan\-te entre deux hypersurfaces
analytiques réelles et strictement pseudoconvexes de $\mC^n$
se prolonge localement biholomorphiquement.
\end{theo}

Ainsi $d$ se prolonge en un biholomorphisme au voisinage des points fixes
de $\phi$ et $\wt{\phi}$.

On applique le principe de prolongement des conjugaisons
pour compléter le remplissage $R$ en une variété complexe sans bord :
grâce au lemme \ref{bassins} on peut prolonger la conjugaison
autour des points fixes à une conjugaison sur la boule fermée en
complétant $R$ par $\BS$.

\begin{prop}
Il existe une surface complexe sans bord $X$ et une application
$i:R\rightarrow X$ biholomorphe sur son image, telles que :
\begin{enumerate}
\item $X\priv i(R)$ est isomorphe à la boule standard ;
\item l'automorphisme $i_\ast\phi$ de $i(R)$ se prolonge en un
      automorphisme de $X$, qu'on note encore $\phi$ ;
\item l'action de $\phi$ sur $X$ est conjuguée, au voisinage de
      $X\priv \interieur{R}$, à l'action de $\wt{\phi}$ sur
      un voisinage de la boule unité fermée de $\CP^2$.
\end{enumerate}
\end{prop}

Remarquons que de la même façon on obtient le résultat complètement
élémentaire suivant.

\begin{prop}
À isomorphisme près,
le seul remplissage convexe de la sphère standard
équivariant relativement à un groupe fermé non compact
est la boule standard fermée.

Le seul remplissage concave fortement équivariant de la sphère standard
est $\CP^2\priv\BS$.
\end{prop}

Le lemme \ref{bassins} 
livre une information cruciale sur $X$.

\begin{lemm}\label{pinceau}
La surface $X$ contient une courbe lisse, irréductible,
rationnelle d'auto-intersection $1$ qui appartient à un pinceau 
préservé par $\phi$
dont au plus un élément est réductible et dont l'unique point-base
est le point fixe attractif de $\phi$ dans $\bord R$.
\end{lemm}

D'après le lemme de Noether, comme $X$ contient en particulier
une courbe rationnelle d'auto-intersection $1$, 
c'est une surface rationnelle.

On note $\vert C\vert$ le pinceau obtenu dans le lemme
\ref{pinceau} et $C$ une de ses courbes génériques.

\subsection{Contraction des courbes exceptionnelles}

On cherche maintenant à contracter les courbes exceptionnelles
de $X$ pour se ramener au cas des surfaces rationnelles
minimales, qu'on comprend bien.

\begin{lemm}
Il existe un entier $k$ tel que $\phi^k$ fixe globalement
chacune des courbes exceptionnelles.
\end{lemm}

\preuve
les diviseurs de $X$ se plongent naturellement dans l'espace
de cohomologie $H^{1,1}_\mR(X)$ sur lequel la forme d'intersection
se prolonge en une forme quadratique symétrique non dégénérée.
Le théorème d'indice de Hodge (voir par exemple \cite{Barth} page
143) dit que cette forme est de signature lorentzienne. En particulier,
puisqu'elle vaut $1$ sur la classe de $C$, elle préserve
son orthogonal $C^\perp$ et sur celui-ci elle 
est définie négative.

Un vecteur de $H^{1,1}_\mR(X)$ se décompose sous la forme
$\lambda C+D$ où $D\in C^\perp$, et la forme définie positive donnée par
\[\lambda^2-D\cdot D\]
est alors préservée par $\phi$.

Considérons l'ensemble ${\mathscr E}$ des courbes exceptionnelles
de $X$. C'est une partie discrète de $H^{1,1}_\mR(X)$
et comme $\phi$ agit par isométrie pour une forme définie positive,
il existe un exposant $k$ tel que l'action de $\phi^k$ sur ${\mathscr E}$
est triviale. Mais une courbe d'auto-itersection n'egative est la seule
représentante complexe de sa classe d'homologie, donc 
chaque courbe exceptionnelle est globalement fixée par
$\phi^k$.
\finpreuve

On veut maintenant contracter les courbes rationnelles 
d'auto-inter\-sec\-tion $-1$ :
comme chacune de ces courbes est préservée, $\phi$ sera rationnellement
conjuguée par cette contraction à un automorphisme d'une surface minimale.

Toutefois, il ne faut pas contracter de courbe qui rencontre $\bord R$
sous peine de singulariser le bord. Le lemme suivant nous permet de contourner
cette difficulté.

\begin{lemm}
Si $E$ est une courbe rationnelle d'auto-intersection $-1$ qui rencontre
$\bord R$, alors il existe une autre courbe rationnelle 
d'auto-intersection $-1$, qui de plus évite $\bord R$.
\end{lemm}

\preuve
on reprend les détails de la démonstration du lemme de Noether
pour déceler cette nouvelle courbe (voir \cite{GriffithsHarris}
page 513).

La courbe $E$, puisqu'elle est
préservée par $\phi$, est nécessairement tangente à $\bord R$ en
un point fixe de $\phi$. De plus la conjugaison à $\CP^2$ doit l'envoyer,
au moins au voisinage de ce point, sur la droite
tangente à $\SSD^3$ en $p_+$ ou $p_-$.

Il s'ensuit que le pinceau $\vert C\vert$ contient bien une 
courbe réductible $C_0$, dont $E$ est l'une des composantes.
Éclatons le point-base de $\vert C\vert$ et notons
$C^\wedge$, $C_0^\wedge$ et $E^\wedge$ les transformées propres. 

Notons
\[C_0^\wedge=\sum_\nu a_\nu C_\nu\]
où les $C_\nu$ sont irréductibles, l'une d'elles est $E^\wedge$ 
et $a_\nu>0$. Maintenant l'auto-intersection
de $C^\wedge$ est $0$ et les $C_\nu$ tout comme $E^\wedge$ sont
disjointes de $C^\wedge$. On en déduit que 
\[C_0^\wedge\cdot C_\nu=C^\wedge\cdot C_\nu=0\]
puis que toutes les composantes $C_\nu$ sont d'auto-intersection
négative. La formule d'adjonction donne, en notant $K$ le fibré
canonique de $X$,
\begin{equation}
\frac{C_0^\wedge\cdot C_0^\wedge+C_0^\wedge\cdot K}{2}+1=0
\end{equation}
et, comme $C_0^\wedge\cdot C_0^\wedge=0$, $C_{\nu_0}\cdot K<0$ pour 
un certain $\nu_0$. Il découle du critère fort de contractibilité
de Castelnuovo--Enriques
que $C_{\nu_0}$ peut être contractée.
Or l'auto-intersection de $E^\wedge$ vaut $-2$, donc $C_{\nu_0}\neq E$
et $C_{\nu_0}$ doit éviter $(\bord R)^\wedge$, en particulier
n'est pas affectée par l'éclatement.
Avant celui-ci il existait donc bien une courbe contractable
évitant $\bord R$.
\finpreuve

Si $X$ n'est pas minimale, on peut donc trouver une courbe à contracter.
On note avec un $\vee$ en exposant l'image d'un objet par la contraction. 
Comme la courbe exceptionnelle contractée est globalement préservée par 
$\phi$, on obtient
 un biholomorphisme $\phi^\vee$ de $X^\vee$ qui fixe 
le point image de la courbe contractée. On peut alors continuer
car $\phi^\vee$ préserve encore une sphère standard, $\bord R^\vee$.
Après un nombre fini d'étapes,
on arrive à une surface rationnelle minimale.

Il suffit maintenant de montrer un dernier lemme pour obtenir le
théorème \ref{theoremplissages}.

\begin{lemm}
Si $X$ est minimale, $X=\CP^2$.
\end{lemm}

\preuve
les surfaces rationnelles minimales sont $\CP^2$ et les surfaces
de Hirzebruch (c'est-à-dire les fibrés en $\CP^1$ sur $\CP^1$)
$\Sigma_n$ pour $n\in\mN$, $n\neq 1$.
Supposons que la surface $X$ est une surface de Hirzebruch.
La classe fondamentale de la courbe $C$ d'auto-intersection $1$
doit alors se décomposer en
\[[C]=a[F]+b[B]\]
où $a$, $b$ sont des entiers, $F$ est une fibre et $B$ est la base.
Or si $X=\Sigma_n$ on a :  
\[F\cdot F=0\qquad F\cdot B=1\qquad B\cdot B=-n\] 
d'où
\[1=C\cdot C=2ab-nb^2=b(2a-nb)\]
On en déduit $b=1$ et $a=(n+1)/2$, d'où
\[C\cdot B=a-nb=\frac{1-n}{2}\]
et donc $n=1$, mais $\Sigma_1$ est justement la seule
surface de Hirzebruch à ne pas être minimale.
\finpreuve

\section{Exemples singuliers}

Dans cette dernière section nous proposons de montrer comment,
en acceptant les sphères à une singularité, on peut construire
des presque-exemples sur $\Sigma_n$.

Considérons une transformation parabolique $\phi$ de $\CP^2$ préservant 
la boule unité, choisie de façon à fixer chaque point d'une droite
projective $L$. Celle-ci est donc tangente à la sphère unité.
Éclatons deux points de $L$, loin du bord et notons $E_1$, $E_2$ les
courbes exceptionnelles obtenues. La transformée propre
$L^\wedge$ est alors d'auto-intersection $-1$, et on peut la contracter.
À travers cette transformation, $\phi$ passe en un automorphisme
de $\CP^1\times\CP^1$ qui préserve un domaine biholomorphe à la boule,
dont le bord est topologiquement une sphère mais possède un point singulier
au point image de la courbe contractée. Ce point est
en effet le seul point d'intersection de la sphère avec les
transformées propre $E_i^\vee$, et si elle était différentiable en
ce point elle devrait leur être tangente à toutes les deux. Mais ces deux
courbes sont chacune une fibre d'un des deux réglages de $\CP^1\times\CP^1$
et sont donc transverses.

\begin{figure}[htb]\begin{center}
\input{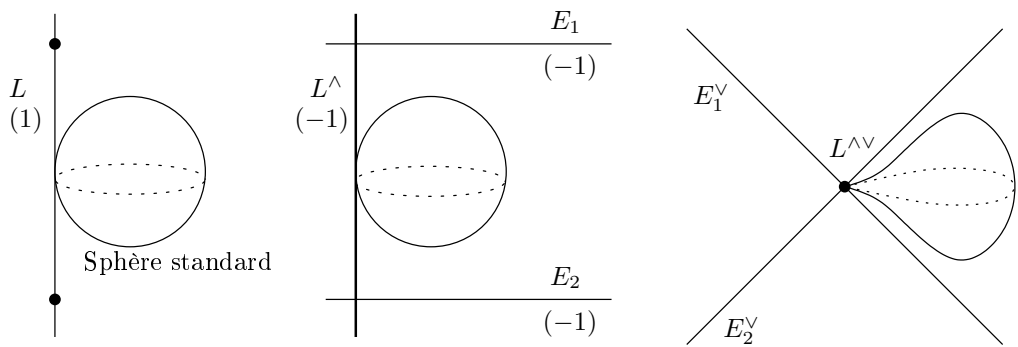}
\caption{Construction d'un exemple singulier dans $\Sigma_0$}
\end{center}\end{figure}

Considérons maintenant une transformation $\phi$ qui possède un 
point fixe hors de la boule unité fermée et éclatons ce point fixe. 
Dans la surface $\Sigma_1$ obtenue la transformée propre $L^\wedge$
de la droite joignant le point éclaté à l'un quelconque des points
fixes de $\phi$ sur la sphère est d'auto-intersection $0$. Le point
d'intersection
entre $L^\wedge$ et la courbe exceptionnelle créée par l'éclatement est
alors fixé par $\phi^\wedge$, conjugué de $\phi$ par l'éclatement.
On peut donc éclater $\Sigma_1$ en ce point puis contracter la transformée
propre de $L^\wedge$ pour obtenir un automorphisme
de $\Sigma_2$ qui préserve un domaine biholomorphe à la boule,
dont le bord est topologiquement une sphère mais possède un point singulier.

\begin{figure}[htb]\begin{center}
\input{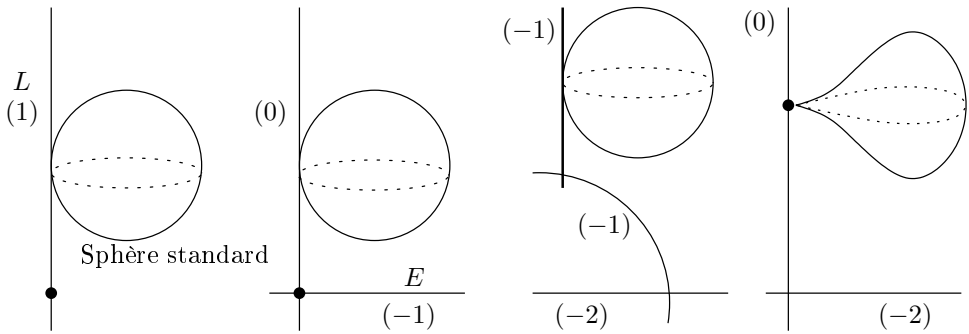}
\caption{Construction d'un exemple singulier dans $\Sigma_2$}
\end{center}\end{figure}

Ensuite on peut passer à $\Sigma_n$ en répétant autant de fois que
nécessaire l'éclatement du point d'intersection entre la
fibre portant le point singulier et la base et la contraction
de la transformée propre de cette fibre.

L'existence de ces exemples, singuliers mais dociles, amène à se
poser la question suivante. Peut-on trouver une surface $X$ possédant
automorphisme qui :
\begin{enumerate}
\item est d'entropie topologique positive et
\item préserve une hypersurface réelle $H$ et
\item induit sur $H$ un groupe non relativement compact
\end{enumerate}
en supposant par exemple que
\begin{enumerate}
\item $H$ est singulière en un ou plusieurs points et strictement
      pseudoconvex là où elle est lisse ou
\item $H$ est lisse partout et strictement pseudoconvexe sauf le long d'une
      sous-variété ?
\end{enumerate}

\appendix
\section*{Annexe : sous-groupes purement elliptiques}
\addcontentsline{toc}{section}{Annexe : sous-groupes purement 
elliptiques}

Dans cette annexe nous nous intéressons à certains sous-groupes de
transformations d'espaces symétriques à courbure négative.
L'objectif principal est de démontrer le théorème \ref{casferme},
utilisé dans la démonstration de \ref{theoremplissages}.

Soit $M$ un espace symétrique à courbure négative et $G$ son groupe
d'isométries. Un élément de $G$ est dit \defini{elliptique} s'il
admet un point fixe dans $M$. Un sous-groupe $F$ de $G$ est dit 
\defini{purement elliptique}\index{Sous-groupe!purement elliptique} 
si tous ses éléments sont elliptiques. On dit que $F$ admet un point
fixe s'il est purement elliptique et s'il existe un point de $M$
fixé simultanément par tous les éléments de $F$.

\begin{theo}\label{casferme}
Un sous-groupe fermé purement elliptique d'isométries 
d'un espace symétrique à courbure négative est compact.
\end{theo}

Dans ce résultat l'hypothèse de fermeture est capitale, comme
le montre le résultat surprenant qui suit.
\begin{theo}[Waterman \cite{Waterman}]
Il existe des sous-groupes d'isométries de $\RH^n$ et de $\mR^{n-1}$
qui sont purement elliptiques et sans point fixe, dès que $n\geqslant 5$.
\end{theo}

Dans le théorème \ref{casferme}, on peut remplacer «~fermé~»
par d'autres hypothèses (voir par exemple \cite{Kulkarni}).
Par exemple, si tous les éléments de $F$ sont d'ordre fini,
alors $F$ est relativement compact. En effet un théorème de
Schur assure qu'un groupe de matrices de torsion est virtuellement
abélien. Alors $F$ contient un sous-groupe distingué $H$, abélien
donc diagonalisable. Comme $H$ est de torsion, il est inclu dans
un tore donc est relativement compact. Mais $F$ est la réunion finie
de classes à gauches toutes homéomorphes à $H$, donc est lui-même
relativement compact.

On va relier la compacité de $F$ et l'existance de point fixe en
commençant par énoncer un résultat naturel dû à Cartan.
\begin{theo}\label{PointFixe}
Un sous-groupe compact de $G$ admet toujours un point fixe.
\end{theo}

On peut facilement exhiber des groupes purement elliptiques et non 
compacts : il suffit par exemple de considérer le groupe engendré 
par une rotation irrationnelle. Toutefois ces groupes sont relativement
compacts, c'est-à-dire que leur adhérence est compacte. Comme les
groupes compacts maximaux de $G$ sont les fixateurs des points de $M$,
on a l'équivalence suivante.

\begin{prop}
Un sous-groupe de $G$ est relativement compact si et seulement
s'il admet un point fixe.
\end{prop}

On démontre maintenant le théorème 
\ref{casferme}, en procédant par étapes.
Le principe est simple : si $F$ est discret 
ou connexe, le résultat découle respectivement du
lemme de Selberg ou d'un théorème
de Montgomery et Zippin. Le cas général se déduit de ses deux
cas extrêmes.

\begin{theo}[Lemme de Selberg]
Un sous-groupe de $\GL{n}{\mR}$ de type fini est virtuellement sans
torsion.
\end{theo}
Par «~virtuellement sans torsion~», on entend qu'il contient un
sous-groupe normal d'indice fini et sans torsion.

\begin{coro}\label{casdiscret}
Si $F$ est purement elliptique et discret, il admet un point fixe dans 
$M$.
\end{coro}

\preuve c'est une généralisation d'une démonstration de \cite{Waterman}.

Tout élément $a$ de $F$ est d'ordre fini : $a$ est elliptique donc
agit comme une rotation sur l'espace tangent de l'un quelconque de
ses points fixes, et comme $F$ est discret les angles de cette rotation
sont commensurables à $\pi$. 

De plus $F$ est dénombrable donc il est
l'union dénombrable de sous-groupes de type fini 
$F_n=\langle a_1,\dots,a_n\rangle$.

Comme $M$ est un espace symétrique, $G$ est un groupe de matrices
et il en est de même pour $F$ et $F_n$. On peut ainsi appliquer
le lemme de Selberg à $F_n$ qui est donc virtuellement sans torsion. 
Mais on a montré que c'est
un groupe de torsion donc il est fini.

D'après le théorème \ref{PointFixe}, il s'ensuit que l'ensemble $V_n$
des points fixes de $F_n$ est non vide. Or $V_n$ est une
sous-variété totalement géodésique et complète de $M$.
Une suite décroissante de telles sous-variétés est stationnaire,
donc $\bigcap V_n\neq\varnothing$ et $F$ admet un point fixe.
\finpreuve

\begin{theo}[Montgomery-Zippin \cite{Montgomery}]\label{MontgomeryZ}
Un groupe de Lie con\-nexe et non compact contient un sous-groupe à
un paramètre fermé et non compact.
\end{theo}

\begin{coro}\label{casconnexe}
Si $F$ est purement elliptique, fermé et connexe, il admet un point
fixe dans $M$.
\end{coro}

\preuve
comme $F$ est purement elliptique, ses sous-groupes fermés à un paramètre
sont tous des cercles. D'après le théorème \ref{MontgomeryZ},
$F$ est donc compact. On conclut par le théorème \ref{PointFixe}.
\finpreuve

On peut maintenant démontrer le théorème \ref{casferme}.
On suppose donc que $F$ est un groupe purement elliptique
et fermé de $G$ et on va montrer qu'il admet un point fixe.

Pour $f\in F$, on note $\text{Fix}(f)$ l'ensemble des points fixes
de $f$.

D'après le corollaire \ref{casconnexe}, la composante neutre
$F_0$ de $F$ est compacte. Alors l'ensemble
\[P_0=\bigcap_{f\in F_0} \text{Fix}(f)\]
des points fixés par tout $F_0$ est non vide et c'est une 
sous-variété complète et totalement
géodésique de $M$. À ce titre, c'est un espace symétrique.

Considérons $g\in F\priv F_0$. Alors on a
\begin{eqnarray}
gP_0 & = & \bigcap_{f\in F_0} g\text{Fix}(f) \nonumber\\
     & = & \bigcap_{f\in F_0} \text{Fix}(gfg^{-1}) \nonumber\\
     & = & P_0 \nonumber
\end{eqnarray}
donc $g$ agit sur $P_0$. De plus, n'importe quel autre élément
de la composante connexe de $g$ dans $F$ s'écrit $gf$ où
$f\in F_0$ donc cette action passe au quotient en une action
du groupe discret $F/F_0$ sur l'espace symétrique $P_0$.
D'après le corollaire \ref{casdiscret}, $F/F_0$ admet un
point fixe dans $P_0$, qui est alors un point fixe de $F$.

\section*{Remerciements}

Ce travail n'aurait jamais vu le jour sans mon directeur de thèse
Abdelghani Zeghib, qui m'a proposé cette question et a su m'encourager
à y répondre. Il ne serait pas ce qu'il est sans les 
explications de Serge Cantat, Jean-Claude Sikorav et Jean-Yves Welschinger.
Merci à tous.

\bibliographystyle{plain}
\bibliography{biblio}

\def\cprime{$'$} \def\dbar{\leavevmode\hbox to 0pt{\hskip.2ex\accent"16\hss}d}
\begin{thebibliography}{10}

\bibitem{Barth}
Wolf~P. Barth, Klaus Hulek, Chris A.~M. Peters, and Antonius Van~de Ven.
\newblock {\em Compact complex surfaces}, volume~4 of {\em Ergebnisse der
  Mathematik und ihrer Grenzgebiete. 3. Folge. A Series of Modern Surveys in
  Mathematics [Results in Mathematics and Related Areas. 3rd Series. A Series
  of Modern Surveys in Mathematics]}.
\newblock Springer-Verlag, Berlin, second edition, 2004.

\bibitem{Eliashberg}
Yakov Eliashberg.
\newblock Filling by holomorphic discs and its applications.
\newblock In {\em Geometry of low-dimensional manifolds, 2 (Durham, 1989)},
  volume 151 of {\em London Math. Soc. Lecture Note Ser.}, pages 45--67.
  Cambridge Univ. Press, Cambridge, 1990.

\bibitem{GriffithsHarris}
Phillip Griffiths and Joseph Harris.
\newblock {\em Principles of algebraic geometry}.
\newblock Wiley Classics Library. John Wiley \& Sons Inc., New York, 1994.
\newblock Reprint of the 1978 original.

\bibitem{SCV}
G.~M. Khenkin, editor.
\newblock {\em Several complex variables. {III}}, volume~9 of {\em
  Encyclopaedia of Mathematical Sciences}.
\newblock Springer-Verlag, Berlin, 1989.
\newblock Geometric function theory, A translation of Sovremennye problemy
  matematiki. Fundamentalnye napravleniya, Tom 9, Akad. Nauk SSSR, Vsesoyuz.
  Inst. Nauchn. i Tekhn. Inform., Moscow, 1986 [ MR0860607 (87m:32011)],
  Translation by J. Peetre, Translation edited by G. M. Khenkin.

\bibitem{Kulkarni}
Ravi~S. Kulkarni.
\newblock Conjugacy classes in {$M(n)$}.
\newblock In {\em Conformal geometry (Bonn, 1985/1986)}, Aspects Math., E12,
  pages 41--64. Vieweg, Braunschweig, 1988.

\bibitem{McMullen}
Curtis~T. McMullen.
\newblock Dynamics on {$K3$} surfaces: {S}alem numbers and {S}iegel disks.
\newblock {\em J. Reine Angew. Math.}, 545:201--233, 2002.

\bibitem{McMullen2}
Curtis~T. McMullen.
\newblock Dynamics on blowups of the projective plane, 2006.
\newblock prépublication.

\bibitem{Montgomery}
Deane Montgomery and Leo Zippin.
\newblock Existence of subgroups isomorphic to the real numbers.
\newblock {\em Ann. of Math. (2)}, 53:298--326, 1951.

\bibitem{Schoen}
R.~Schoen.
\newblock On the conformal and {CR} automorphism groups.
\newblock {\em Geom. Funct. Anal.}, 5(2):464--481, 1995.

\bibitem{Waterman}
P.~L. Waterman.
\newblock Purely elliptic {M}\"obius groups.
\newblock In {\em Holomorphic functions and moduli, Vol.\ II (Berkeley, CA,
  1986)}, volume~11 of {\em Math. Sci. Res. Inst. Publ.}, pages 173--178.
  Springer, New York, 1988.

\end{thebibliography}

\end{document}